\documentclass[12pt,a4paper]{article}
\usepackage[round]{natbib}
\usepackage{latexsym,amsfonts,amssymb,amsmath,amsthm,amsxtra}

\begin{document}

\vspace*{0.9cm}
\begin{center}
{\bf MODERATE DEVIATION THEOREM \\FOR THE NEYMAN-PEARSON STATISTIC\\ IN TESTING UNIFORMITY}

\vspace{1cm}
{\bf Tadeusz Inglot}\\

\vspace{5mm}
{\it Faculty of Pure and Applied Mathematics\\
Wroc{\l}aw University of Science and Technology}
\end{center}

\vspace{0.3cm}
\begin{quotation}
{\small{\bf Abstract.} We show that for local alternatives to uniformity which are determined by a sequence of square integrable densities the moderate deviation (MD) theorem for the corresponding Neyman-Pearson statistic does not hold in the full range for all unbounded densities. We give a sufficient condition under which MD theorem holds. The proof is based on Mogulskii's inequality.\\

\noindent{\it Key words and phrases: testing for uniformity, local alternatives, Neyman-Pearson statistic, moderate deviations, square integrable density, Mogulskii's inequality.\\
MSC Subject Classification:62G10, 60F10, 62G20.}}

\end{quotation}

\vspace{0.4cm}
\noindent{\bf 1. Introduction}\\

The intermediate approach to tests' comparison was initiated by \citet*{r14} and developed by \citet*{r12}, \citet*{r8, r10}, \citet*{r7}, \citet*{r3, r5, r6}, \citet*{r11}, among others. Similarly as for the Bahadur efficiency, the intermediate efficiency is calculated as a limit of the ratio between two  slopes. The intermediate slope is determined by an index of moderate deviations under the null hypothesis and a scalling factor resulting from a kind of weak law of large numbers under the sequence of alternatives. By an index of moderate deviations (MD) for a generic statistic $T_n$ we mean the limit
$$-\lim_{n\to\infty}\frac{1}{nx_n^2}\log Pr(T_n\geqslant \sqrt{n}x_n)=c,\eqno (1)$$
provided it exists and is positive, where $Pr$ represents a null distribution while $x_n$ are positive, $x_n\to 0$ and $nx_n^2\to\infty$ as $n\to\infty$. The relation (1) we shall call MD theorem for $T_n$.

The Neyman-Pearson test seems to be the most natural procedure to which other tests could be compared. MD theorem for the Neyman-Pearson statistic in the full range i.e. for all $x_n\to 0$ such that $nx_n^2\to\infty$ as $n\to\infty$ is one of sufficient conditions to make it possible (cf. \citealp{r11}, \citealp{r2}).

In the present paper we study this last question in the classical case of testing for uniformity.

Let $X_1,X_2,..., X_n$ be a sample from a distribution $P$ on the interval $[0,1]$. Consider testing
$$H_0: P=P_0,$$
where $P_0$ is the uniform distribution over $[0,1]$. Let $P_n$ be a sequence of local alternatives, convergent to $P_0$, given by densities $p_n(t)=1+\vartheta_na(t)$, where $\vartheta_n\to 0$ while $a\in L_2(0,1)$ is fixed and satisfies
$$\int_0^1 a(t)dt=0,\;\;\int_0^1 a^2(t)dt=1.\eqno (2)$$
The normalized Neyman-Pearson statistic for testing $H_0$ against the alternative with density $p_n$ has the form
$$ V_n=\frac{1}{\sqrt{n}\sigma_{0n}}\sum_{i=1}^n (\log(1+\vartheta_n a(X_i))-e_{0n}), \eqno (3)$$
where $e_{0n}=\int_0^1 \log(1+\vartheta_n a(t))dt$ and $\sigma_{0n}^2=\int_0^1\log^2(1+\vartheta_na(t))dt-e_{0n}^2$ are normalizing sequences. 

In the paper by \citet*{r8} it was proved that for $V_n$ with $a$ bounded (1) holds in the full range of sequences $x_n$ (Theorem 1, below). In many typical goodness of fit testing problems like e.g. testing in the Gaussian shift or the Gaussian scale families the transformation onto $(0,1)$ leads to unbounded or even not square integrable functions $a$ (see e.g. \citealp{r2}, section 8).
Our main result (Theorem 2 and Corollary) gives sufficient conditions on $x_n$ under which (1) holds for $V_n$. We also show (Theorem 3) that (1) does not hold for $V_n$ in the full range of $x_n$ at least for some unbounded functions which can belong to $L_q(0,1)$ with arbitrary $q>2$. All proofs are sent to Section 3.

Throughout the rest of the paper we assume that $H_0$ is true i.e. that $X_i$ are uniformly distributed over $[0,1]$. Also by $P_0^n$ we denote $n$-fold product of $P_0$ and by $E_0$ and Var$_0$ an expectation and a variance calculated under $P_0$ or $P_0^n$.\\

\vspace{2mm}
\noindent{\bf 2. Moderate dviations for $V_n$}\\

We start with asymptotic formulae for normalizing sequences $e_{0n},\;\sigma_{0n}$ in (3) which will be exploited in the sequel.\\

{\bf Proposition 1.} {\it If $a\in L_2(0,1)$ then} 
$$e_{0n}=-\frac{\vartheta_n^2}{2}(1+o(1)) \eqno (4)$$
and
$$\sigma_{0n}=\vartheta_n(1+o(1)). \eqno (5)$$

\vspace{2mm}
Now, assume that $a$ in (3) is bounded. Theorem 1, below, recalls the MD theorem for $V_n$ for bounded $a$ obtained in \citet*{r8}. In that paper it was proved using MD result for triangular arrays of independent random variables from the unpublished paper by \citet*{r1}. In Section 3 we reprove this theorem by reducing  to the classical MD theorem for i.i.d bounded random variables.\\

{\bf Theorem 1.} {\it Suppose $|a|\leqslant M$ for some $M\geqslant 1$. Then for every positive $x_n$ such that $x_n\to 0$ and $nx_n^2\to\infty$ we have}
$$ -\lim_{n\to\infty}\frac{1}{nx_n^2}\log P_0^n(V_n\geqslant \sqrt{n}x_n)=\frac{1}{2}.\eqno (6)$$

\vspace{2mm}
Next, suppose that $a\in L_2(0,1)$ in (3) is unbounded. Under this assumption we are able to get (6) for $x_n$ satisfying some additional restriction. The proof goes along the same line of argument as that for the classical MD theorem for i.i.d. random variables based on a version of Mogulskii's inequality \citep{r13} proposed in \citet*{r4}. Therefore in the Appendix we provide the proof of this classical theorem (Theorem 4) to show that indeed large parts of the proof of Theorem 2 are simply rewriting those of Theorem 4.\\

{\bf Theorem 2.} {\it Suppose $a\in L_2(0,1)$ is unbounded and $\vartheta_n\to 0$ is such that $n\vartheta_n^2\to\infty$. \\
(i)\quad For any $\delta>0$ and every positive $x_n$ such that $x_n\leqslant (1-\delta)\sigma_{0n}$ and $nx_n^2\to\infty$ we have
$$ -\limsup_{n\to\infty}\frac{1}{nx_n^2}\log P_0^n(V_n\geqslant \sqrt{n}x_n)\geqslant \frac{1}{2};$$
(ii)\quad for any $\delta>0$ and every positive $x_n$ such that $x_n\leqslant \frac{1}{3}(1-\delta)\sigma_{0n}$ and $nx_n^2\to\infty$ we have}
$$ -\liminf_{n\to\infty}\frac{1}{nx_n^2}\log P_0^n(V_n\geqslant \sqrt{n}x_n)\leqslant \frac{1}{2}.$$

\vspace{3mm}
Theorem 2 and (5) immediately imply the following corollary.\\

{\bf Corollary.} {\it Suppose $a\in L_2(0,1)$ is unbounded and $\vartheta_n\to 0$ is such that $n\vartheta_n^2\to\infty$. 
Then for every positive $x_n$ such that $\limsup_{n\to\infty}(x_n/\vartheta_n)<1/3$ and $nx_n^2\to\infty$ the relation (6) holds.}\\

Denote random variables $Y_{ni}=(\log(1+\vartheta_na(X_i))-e_{0n})/\sigma_{0n},\;i=1,...,n.$ Then $E_0Y_{ni}=0$, Var$_0Y_{ni}=1$ and $\varphi_n(\lambda)=E_0\exp\{\lambda Y_{ni}\}=e^{-\lambda e_{0n}/\sigma_{0n}}E(1+\vartheta_na(X_i))^{\lambda/\sigma_{0n}}<\infty$ for $\lambda\leqslant 2\sigma_{0n}$. \\

{\bf Remark.} If $a\notin L_q(0,1)$ for some $q>2$ then $\varphi_n(\lambda)=\infty$ for $\lambda\geqslant q\sigma_{0n}$. Therefore for $a\in L_2(0,1)$ not belonging to $L_q(0,1)$ for all $q>2$ the moment generating function $\varphi_n(\lambda)$ does not exists when $\lambda/\vartheta_n$ is sufficiently large. This suggests that Theorem 2 and Corollary cannot be essentially strenghtened and the condition $\limsup_{n\to\infty}x_n/\vartheta_n<\infty$ seems to be necessary for (6). The next theorem partially confirms such a conjecture.\\

Consider unbounded square integrable functions satisfying (2) of the form
$$ a_r(t)=\frac{\sqrt{1-2r}}{r}\left(\frac{1-r}{t^r}-1\right),\;\;r\in \left(0,\frac{1}{2}\right),$$
corresponding sequences of local alternatives and the Neyman-Pearson statistics (3).\\ 

{\bf Theorem 3.} {\it Suppose $V_n$ is the Neyman-Pearson statistic (3) applied to the function $a_r$ for some $r\in (0,1/2)$ and $\vartheta_n\to 0$ with $n\vartheta_n^2\to\infty$. If positive $x_n$ fulfill the following condition
$$ \mbox{for some}\;\;q<r\;\;\mbox{it holds}\;\;\frac{x_n}{\vartheta_n^q}\to \infty\;\;\mbox{and}\;\;x_n^{(r-q)/q}\log\vartheta_n\to 0$$
then}
$$\lim_{n\to\infty}\frac{1}{nx_n^2}\log P_0^n(V_n\geqslant \sqrt{n}x_n)=0.\eqno(7)$$

\vspace{2mm}
Theorem 3 shows that in every space $L_q(0,1),\;q>2,$ there are functions satisfying (2) such that (6) does not hold for all $x_n\to 0$ such that $nx_n^2\to\infty$. This means that Theorem 1 can not be extended to the class of all square integrable functions $a$. 

Theorem 2 applied to the function $a_r$ and Theorem 3 do not cover a wide range of sequences $x_n$ for which validity of (6) for this particular $a_r$ remains undecided.\\

\noindent{\bf 3. Proofs}\\

\noindent{\bf Proof of Proposition 1.} Let $\varepsilon\in (0,1)$ be arbitrary. Then the inequality 
$$  y-\frac{3-\varepsilon}{6(1-\varepsilon)}y^2\leqslant \log(1+y)\leqslant y-\frac{3-2\varepsilon}{6}y^2 \eqno (8)$$
holds on $[-\varepsilon,\varepsilon]$. From Markov's inequality we have $P_0(a^2(X_1)>\varepsilon^2/\vartheta_n^2)\leqslant \vartheta_n^2/\varepsilon^2$. Hence and from the Cauchy-Schwarz inequality we obtain for $n$ sufficiently large (i.e. such that $\vartheta_n<\varepsilon$)
$$ \int_{\vartheta_na>\varepsilon} a(t) dt\leqslant \sqrt{\int_{\vartheta_na>\varepsilon} a^2(t)dt}\sqrt{P_0(\{\vartheta_na(X_1)>\varepsilon\})}\leqslant \frac{\vartheta_n}{\varepsilon }o(1).\eqno (9)$$
So, from (2), (8) and (9) we get
$$e_{0n}=\int_0^1 \log(1+\vartheta_na(t))dt\geqslant \int_{\vartheta_na\leqslant\varepsilon}\log(1+\vartheta_na(t))dt\hspace*{2cm}$$
$$\geqslant \int_{\vartheta_na\leqslant\varepsilon}\vartheta_n a(t) dt-\frac{3-\varepsilon}{6(1-\varepsilon)}\vartheta_n^2\int_{\vartheta_na\leqslant\varepsilon}a^2(t)dt\hspace*{1.7cm}$$
$$\hspace*{1cm}\geqslant-\vartheta_n\int_{\vartheta_na>\varepsilon}a(t) dt-\frac{3-\varepsilon}{6(1-\varepsilon)}\vartheta_n^2
\geqslant -\frac{\vartheta_n^2}{\varepsilon}o(1)-\frac{3-\varepsilon}{6(1-\varepsilon)}\vartheta_n^2.$$
Similarly, from (2) and (8) we get
$$e_{0n}=\int_{\vartheta_na\leqslant\varepsilon}\log(1+\vartheta_na(t))dt+\int_{\vartheta_na>\varepsilon}\log(1+\vartheta_na(t))dt\hspace*{4cm}$$
$$\hspace*{0.3cm}\leqslant-\vartheta_n\int_{\vartheta_na>\varepsilon}a(t) dt-\frac{3-2\varepsilon}{6}\vartheta_n^2\int_{\vartheta_na\leqslant\varepsilon}a^2(t)dt+\int_{\vartheta_na>\varepsilon}\log(1+\vartheta_na(t))dt$$
$$\leqslant -\vartheta_n\int_{\vartheta_na>\varepsilon}a(t) dt-\frac{3-2\varepsilon}{6}\vartheta_n^2(1+o(1))+\vartheta_n\int_{\vartheta_na>\varepsilon}a(t) dt=-\frac{3-2\varepsilon}{6}\vartheta_n^2(1+o(1)).$$
Hence for arbitrary $\varepsilon\in (0,1)$ we have
$$-\frac{3-\varepsilon}{6(1-\varepsilon)}\leqslant\liminf_n \frac{e_{0n}}{\vartheta_n^2}\leqslant \limsup_n \frac{e_{0n}}{\vartheta_n^2}\leq -\frac{3-2\varepsilon}{6}.$$
Since $\varepsilon$ is arbitrary (4) follows. 

In the same way we show (5) (cf. Proposition 3 in \citealt{r6}).\hfill $\Box$\\

\noindent{\bf Proof of Theorem 1.} On $(-1,\infty)$ define a function $h(y)=2\frac{y-\log(1+y)}{y^2}$ with $h(0)=1$. The function $h(y)$ is of class $C^{\infty}$, positive and decreasing on $(-1,\infty)$ and analytic on $(-1,1)$. Since $\log(1+y)=y-\frac{y^2}{2}h(y)$ then from (2) we get
$$ e_{0n}=E_0\log(1+\vartheta_n a(X_1))=-\frac{\vartheta_n^2}{2}E_0a^2(X_1)h(\vartheta_na(X_1))=-\frac{\vartheta_n^2}{2}\mu_n,$$
where $\mu_n=1+o(1)$ from (4) (or from Lebesgue's Dominated Convergence Theorem). This implies
$$P_0^n(V_n\geqslant \sqrt{n}x_n)\hspace*{4cm}$$
$$= P_0^n\left(\left[\frac{1}{\sqrt{n}}\sum_{i=1}^n a(X_i)-\frac{\vartheta_n}{2\sqrt{n}}\sum_{i=1}^n\left(a^2(X_i)h(\vartheta_na(X_i))-\mu_n\right)\right]\geqslant\sqrt{n}x_n\frac{\sigma_{0n}}{\vartheta_n}\right).$$
Since $a\geqslant -1$ a.s. then for $n$ sufficiently large random variables $a^2(X_i)h(\vartheta_na(X_i))$ are bounded by $3M^2/2$. Moreover, for $\tau^2_n=\mbox{Var}_0 a^2(X_i)h(\vartheta_na(X_i))$ we have $\tau^2_n\to\int_0^1a^4(t)dt-1$ from (2) and Lebesgue's Dominated Convergence Theorem. Denote
$$ D_n=\left\{\left|\frac{1}{\sqrt{n}}\sum_{i=1}^n \left(a^2(X_i)h(\vartheta_na(X_i))-\mu_n\right)\right|<2\tau_n\sqrt{n}x_n\frac{\sigma_{0n}}{\vartheta_n}\right\}.$$
Then from the classical Bernstein inequality we get
$$ P_0^n(D_n^c)\leqslant 2\exp\left\{-2nx_n^2\frac{\sigma_{0n}^2}{\vartheta_n^2}\frac{1}{1+2M^2x_n\sigma_{0n}/\vartheta_n\tau_n}\right\}=2\exp\{-2nx_n^2(1+o(1))\},$$
where $A^c$ denotes the complement of a set $A$. Hence and denoting $F_n=\{V_n\geqslant \sqrt{n}x_n\}$ we obtain 
$$P_0^n(V_n\geqslant \sqrt{n}x_n)\geqslant P_0^n(F_n\cap D_n)\geqslant P_0^n\left(\left\{\frac{1}{\sqrt{n}}\sum_{i=1}^n a(X_i)\geqslant (1+\tau_n\vartheta_n)\sqrt{n}x_n\frac{\sigma_{0n}}{\vartheta_n}\right\}\cap D_n\right)$$
$$\geqslant P_0^n\left(\frac{1}{\sqrt{n}}\sum_{i=1}^n a(X_i)\geqslant (1+\tau_n\vartheta_n)\sqrt{n}x_n\frac{\sigma_{0n}}{\vartheta_n}\right)-P_0^n(D_n^c)$$
As $\sigma_{0n}/\vartheta_n=1+o(1)$ by (5) then from the classical MD theorem (Theorem 4 in the Appendix) applied to the sequence $a(X_i)$ of bounded random variables the last expression can be estimated from below by
$$\exp\left\{-\frac{nx_n^2}{2}(1+o(1))\right\}-2\exp\{-2nx_n^2(1+o(1))\}. \eqno (10)$$
Similarly
$$P_0^n(V_n\geqslant\sqrt{n}x_n)\leqslant P_0^n(F_n\cap D_n) +P_0^n(D_n^c)$$
$$\leqslant P_0^n\left(\left\{\frac{1}{\sqrt{n}}\sum_{i=1}^n a(X_i)\geqslant (1-\tau_n\vartheta_n)\sqrt{n}x_n\frac{\sigma_{0n}}{\vartheta_n}\right\}\cap D_n\right)+P_0^n(D_n^c)$$
$$\leqslant P_0^n\left(\frac{1}{\sqrt{n}}\sum_{i=1}^n a(X_i)\geqslant (1-\tau_n\vartheta_n)\sqrt{n}x_n\frac{\sigma_{0n}}{\vartheta_n}\right)+P_0^n(D_n^c)\hspace*{1.7cm}$$
$$\leqslant \exp\left\{-\frac{nx_n^2}{2}(1+o(1))\right\}+2\exp\{-2nx_n^2(1+o(1))\}.\hspace*{2.2cm} \eqno (11)$$
From (10) and (11) the relation (6) immediately follows. \hfill $\Box$\\

\noindent{\bf Proof of Theorem 2.} The function $\log^{k}(y)/y^2,\;y\geqslant 1,\; k\geqslant 3,$ is bounded from above by $(k/2)^ke^{-k}$ while 
$$v_k(y)=\frac{\log^k(1+y)}{y^2}=\frac{\log^k(1+y)}{(1+y)^2}\frac{(1+y)^2}{y^2},\;y>0,\;k\geqslant 3,$$
is increasing on the interval $(0,1)$. Therefore $v_k(y)$ is bounded from above by $4(k/2)^ke^{-k}$. Hence and from (4)  for $n$ sufficiently large
$$E_0|Y_{ni}|^k=\int_{a\leqslant 1/\sqrt{\vartheta_n}}\frac{|\log(1+\vartheta_n a(t))-e_{0n}|^k}{\sigma_{0n}^k}dt+\int_{a> 1/\sqrt{\vartheta_n}}\left(\frac{\log(1+\vartheta_na(t))-e_{0n}}{\sigma_{0n}}\right)^kdt$$
$$\leqslant \left(\frac{\log(1+\sqrt{\vartheta_n)}-e_{0n}}{\sigma_{0n}}\right)^{k-2} +4k^ke^{-k}\frac{\vartheta_n^2}{\sigma_{0n}^k}\int_{a> 1/\sqrt{\vartheta_n}}a^2(t)dt$$
and from Stirling's formula for $k\geqslant 3$ and $n$ sufficiently large
$$\frac{E_0|Y_{ni}|^k}{k!}\leqslant \frac{1}{6}\left(\frac{\sqrt{\vartheta_n}-e_{0n}}{\sigma_{0n}}\right)^{k-2}+\frac{4}{\sqrt{2\pi k}}\frac{\vartheta_n^2}{\sigma_{0n}^k}\int_{a> 1/\sqrt{\vartheta_n}}a^2(t)dt\leqslant \frac{\vartheta_n^2}{\sigma_{0n}^k}\omega_n, \eqno (12)$$ 
where $\omega_n=\sigma_{0n}^2(\sqrt{\vartheta_n}-e_{0n})/6\vartheta_n^2+\int_{a> 1/\sqrt{\vartheta_n}}a^2dt=o(1)$.

The function $\varphi_n(\lambda)$ is analytic on the interval $[0,2\sigma_{0n}]$ and $\varphi_n(\lambda)=1+\frac{\lambda^2}{2}\psi_n(\lambda)$, where $\psi_n(\lambda)=1+2\sum_{k=3}^{\infty}\frac{EY_{ni}^{k}}{k!}\lambda^{k-2}$. By (12) we have for $n$ sufficiently large
$$|\psi_n(\lambda)-1|\leqslant 2\sum_{k=3}^{\infty}\frac{\vartheta_n^2}{\sigma_{0n}^k}\omega_n\lambda^{k-2} \eqno (13)$$
and
$$|\psi'_n(\lambda)|\leqslant 2\sum_{k=3}^{\infty}\frac{\vartheta_n^2}{\sigma_{0n}^k}\omega_n(k-2)\lambda^{k-3}. 
\eqno (14)$$

\vspace{3mm}
{\it Proof of (i) (upper estimate).} By Markov's inequality we have for $\lambda\in(0,2\sigma_{0n})$
$$ P_0^n(V_n\geqslant\sqrt{n}x_n)=P_0^n\left(\frac{1}{\sqrt{n}}\sum_{i=1}^n Y_{ni}\geqslant \sqrt{n}x_n\right)$$
$$=P_0^n\left(\prod_{i=1}^n e^{\lambda Y_{ni}}\geqslant e^{n\lambda x_n}\right)\leqslant e^{-n\lambda x_n}\varphi_n^n(\lambda).$$
Putting $\lambda=x_n$ the right hand side takes the form $ e^{-nx_n^2}\varphi_n^n(x_n).$ Since $x_n\leqslant (1-\delta)\sigma_{0n}$ then (13) implies for $n$ sufficiently large
$$\varphi_n(x_n)\leqslant 1+\frac{x_n^2}{2}+x_n^2\frac{\vartheta_n^2}{\sigma_{0n}^2}\omega_n\sum_{k=3}^{\infty}(1-\delta)^{k-2}\leqslant 1+\frac{x_n^2}{2}\left(1+\frac{2\vartheta_n^2\omega_n}{\delta\sigma_{0n}^2}\right)$$
and in consequence
$$\frac{1}{nx_n^2}\log P_0^n(V_n\geqslant \sqrt{n}x_n)\leqslant -\frac{1}{2}+\frac{1}{\delta}\frac{\vartheta_n^2\omega_n}{\sigma_{0n}^2},\eqno(15)$$ 
which completes the proof of (i).

{\it Proof of (ii) (lower estimate).} Denote by $P_n$ the distribution of $Y_{ni}$ and let $Q_{n\lambda}\ll P_n$ be such that $\frac{dQ_{n\lambda}}{dP_n}(y)=e^{\lambda y}/\varphi_n(\lambda)$. Then
$$m_n(\lambda)=\int y dQ_{n\lambda}=\frac{1}{\varphi_n(\lambda)}\int ye^{\lambda y}dP_n(y)=\frac{\varphi'_n(\lambda)}{\varphi_n(\lambda)}$$
and the entropy distance (Kullback -Leibler) of $Q_{n\lambda}$ from $P_n$ is equal to
$$D(Q_{n\lambda}||P_n)=\int\frac{1}{\varphi_n(\lambda)}e^{\lambda y}(\lambda y-\log \varphi_n(\lambda))dP_n(y)=\lambda\frac{\varphi'_n(\lambda)}{\varphi_n(\lambda)}-\log\varphi_n(\lambda).$$
For $n\geqslant 1$ and $\varepsilon\in (0,\delta)$ let $\lambda_n>0$ be such that $m_n(\lambda_n)=(1+\varepsilon)x_n$. Observe that $\lambda_n$ is correctly defined and 
$$\lambda_n<5\sigma_{0n}/6.\eqno(16)$$ 
Indeed, the inequality $(1+y)^{5/6}\geqslant 1+5y/6-y^2/9$, which holds on $[-1/2,\infty)$, and (2) give $\varphi_n(5\sigma_{0n}/6)=e^{-5e_{0n}/6}(1-\vartheta_n^2/9)$. This, convexity of $\varphi_n(\lambda)$, 
the assumption $x_n\leqslant (1-\delta)\sigma_{0n}/3$, (4) and (5) imply for $n$ sufficiently large
$$m_n\left(\frac{5}{6}\sigma_{0n}\right)=\frac{\varphi'_n(\frac{5}{6}\sigma_{0n})}{\varphi_n(\frac{5}{6}\sigma_{0n})}\geqslant \frac{\varphi_n(5\sigma_{0n}/6)-\varphi_n(0)}{(5\sigma_{0n}/6)\varphi_n(5\sigma_{0n}/6)}\geqslant \frac{1}{5\sigma_{0n}/6}-\frac{e^{5e_{0n}/6}}{(5\sigma_{0n}/6)(1-\vartheta_n^2/9)}$$
$$\geqslant\frac{\sigma_{0n}}{3}>(1-\delta^2)\frac{\sigma_{0n}}{3}\geqslant (1+\delta)x_n>(1+\varepsilon)x_n=m_n(\lambda_n),$$
which implies (16) (the function $m_n(\lambda)$ is increasing since $\log\varphi_n(\lambda)$ is strictly convex).
Inserting $\lambda=\lambda_n$ to (13) and (14) and using (16) we get for $n$ sufficiently large
$$|\psi_n(\lambda_n)-1|\leqslant 10\frac{\vartheta_n^2}{\sigma_{0n}^2}\omega_n \;\;\mbox{and}\;\;|\lambda_n\psi'_n(\lambda_n)|\leqslant 60\frac{\vartheta_n^2}{\sigma_{0n}^2}\omega_n.$$
Hence for $n$ sufficiently large (i.e. such that $|\psi_n(\lambda_n)-1|<\varepsilon/8,$ $|\lambda_n\psi'_n(\lambda_n)|<\varepsilon/4,$  $\lambda_n^2\psi_n(\lambda_n)<\varepsilon/4$) we obtain
$$ (1+\varepsilon)x_n=m_n(\lambda_n)=\frac{\varphi'_n(\lambda_n)}{\varphi_n(\lambda_n)}=\frac{\lambda_n\psi_n(\lambda_n)+\lambda_n^2\psi'_n(\lambda_n)/2}{1+\lambda_n^2\psi_n(\lambda_n)/2}\leqslant \lambda_n(1+\varepsilon/4)\leqslant \lambda_n(1+\varepsilon)$$
and similarly
$$ (1+\varepsilon)x_n\geqslant \lambda_n\frac{1-\varepsilon/4}{1+\varepsilon/8}$$
which gives 
$$x_n\leq\lambda_n\leq x_n \frac{(1+\varepsilon)(1+\varepsilon/8)}{1-\varepsilon/4}\leq (1+2\varepsilon)x_n.\eqno (17)$$

For $\lambda_n$ defined above we have
$$ D(Q_{n\lambda_n}||P_n)=(1+\varepsilon)\lambda_nx_n-\log\varphi_n(\lambda_n).$$
Now, we apply the following version of Mogulskii's inequality (\citealt{r13}, cf. Corollary 1 in \citealt{r4}).\\

{\bf Theorem A.} {\it Let $Q\ll P$ and $\xi_1,...,\xi_n$ be i.i.d. random variables with distribution $P$ and $\eta_1,...,\eta_n$ i.i.d. random variables with distribution $Q$. Then for every Borel set $A$, any $M\in \mathbb{R}$ and any $n\geqslant 1$ it holds
$$Pr\left(\frac{\xi_1+...+\xi_n}{n}\in A\right)(1-e^{-M})+e^{-M}\geqslant \exp\{-nD(Q||P)-Mp_n\},\eqno (18)$$
where $p_n=Pr(\eta_1+...+\eta_n\in nA^c)$.}\\

In Theorem A we set $P=P_n,\;Q=Q_{n\lambda_n},\;A=[x_n,\infty),\; M=2nx_n^2$. Observe that the variance of $Q_{n\lambda_n}$ is equal to $\rho_n^2=\varphi''_n(\lambda_n)/\varphi_n(\lambda_n)-m_n^2(\lambda_n)\to 1$ since, similarly as above, from (16) we obtain $|\varphi''_n(\lambda_n)-1|\leq 70\frac{\vartheta_n^2}{\sigma_{0n}^2}\omega_n$. Hence for $n$ sufficiently large, by the assumption $nx_n^2\to\infty$ and from Cantelli's inequality we obtain
$$ p_n=Pr(\eta_1+...+\eta_n< nx_n)=Pr\left(\sum_{i=1}^n(\eta_i-m_n(\lambda_n))< -\varepsilon nx_n\right)\leqslant \frac{n\rho_n^2}{n\rho_n^2+\varepsilon^2n^2x_n^2}\to 0$$
and in consequence from (17) and (18) for $n$ sufficiently large
$$ P_0^n(Y_{n1}+...+Y_{nn}\geqslant nx_n)(1-e^{-2nx_n^2})$$
$$\geqslant \exp\{-(1+\varepsilon)n\lambda_nx_n+n\log(1+\lambda_n^2\psi_n(\lambda_n)/2)-2nx_n^2p_n\}-e^{-2nx_n^2}$$
$$\geqslant \exp\left\{-\frac{1+3\varepsilon}{2}n\lambda_nx_n-2nx_n^2p_n\right\}-e^{-2nx_n^2}\geqslant \exp\left\{(-\frac{1}{2}-\frac{7}{2}\varepsilon)nx_n^2-2nx_n^2p_n\right\}-e^{-2nx_n^2}.$$
Logarithming both sides and dividing by $nx_n^2$ we get
$$\frac{1}{nx_n^2}\log P_0^n(V_n\geqslant \sqrt{n}x_n)\geqslant-\frac{1}{2}-\frac{7}{2}\varepsilon +o(1)$$
which, due to arbitrariness of $\varepsilon$, 
ends the proof of (ii) as well as that of Theorem 2. \hfill $\Box$\\

\noindent{\bf Proof of Theorem 3.} Let $\Gamma_n$ be the distribution on $(0,1)$ with the density
$$g_n(t)=1+\frac{x_n^{(r+q)/q}}{\vartheta_n}{\bf 1}_{(\vartheta_n,2\vartheta_n)}(t)-x_n^{(r+q)/2q}{\bf 1}_{(1-x_n^{(r+q)/2q},1)}(t),$$ 
where ${\bf 1}_A(t)$ denotes the indicator of a set $A$. An elementary calculation gives \linebreak[4]$D(\Gamma_n||P_0)= x_n^{(r+q)/q}\log(x_n^{(r+q)/q}/\vartheta_n)(1+o(1))$. 

Similarly as previously denote $Y_{ni}=(\log(1+\vartheta_na_r(X_i))-e_{0n})/\sigma_{0n}$, $i=1,...,n,$ their distributions by $P_{nr}$ when $X_i$ are uniformly distributed over $(0,1)$, or by $Q_{nr}$ when $X_i$ have the distribution $\Gamma_n$. Since $Y_{ni}$ are bijective (decreasing) functions of $X_i$ then $D(Q_{nr}||P_{nr})=D(\Gamma_n||P_0)=x_n^{(r+q)/q}\log(x_n^{(r+q)/q}/\vartheta_n)(1+o(1))$. 

As $a_r(t)<0$ for $t>(1-r)^{1/r}$ then for $n$ sufficiently large we have
$$E_{\Gamma_n}Y_{n1}\geqslant \frac{x_n^{(r+q)/q}}{\sigma_{0n}\vartheta_n}\int_{\vartheta_n}^{2\vartheta_n}\log(1+\vartheta_na_r(t))dt-\frac{x_n^{(r+q)/2q}}{\sigma_{0n}}\int_{1-x_n^{(r+q)/2q}}^1\log(1+\vartheta_na_r(t))dt$$
$$\geqslant \frac{x_n^{(r+q)/q}}{\sigma_{0n}}\log(1+\vartheta_na_r(2\vartheta_n))\geqslant \frac{\sqrt{1-2r}}{2}\frac{x_n^{(r+q)/q}}{\vartheta_n^r}=\frac{\sqrt{1-2r}}{2}x_n\left(\frac{x_n}{\vartheta_n^q}\right)^{r/q}=\kappa_n\eqno (20)$$
and
$$E_{\Gamma_n}Y_{n1}^2\leqslant \frac{1}{\sigma_{0n}^2}\left(\sigma_{0n}^2+\frac{x_n^{(r+q)/q}}{\vartheta_n}\int_{\vartheta_n}^{2\vartheta_n}\left(\log(1+\vartheta_na_r(t))dt-e_{0n}\right)^2\right)$$
$$\leqslant 1+\frac{x_n^{(r+q)/q}}{\sigma_{0n}^2}\left(\log(1+\vartheta_na_r(\vartheta_n))-e_{0n}\right)^2\leqslant \frac{1}{r^2}\frac{x_n^{(r+q)/q}}{\vartheta_n^{2r}}(1+o(1)).$$
In Mogulskii's inequality set $P=P_{nr},\;Q=Q_{nr},\; M=nx_n^2, \;A=[x_n,\infty)$. From the assumption on $x_n$ and (20) it follows $x_n-\kappa_n<0$ for $n$ sufficiently large. So, by Cantelli's inequality for $n$ sufficiently large
$$p_n=Pr(\eta_1+...+\eta_n<nx_n)\leqslant Pr\left(\sum_{i=1}^n(\eta_i-E_{\Gamma_n}Y_{ni})< n(x_n-\kappa_n)\right)$$
$$\leqslant\frac{nE_{\Gamma_n}Y_{n1}^2}{nE_{\Gamma_n}Y_{n1}^2+n^2(\kappa_n-x_n)^2}\leqslant\frac{x_n^{(r+q)/q}(1+o(1))}{x_n^{(r+q)/q}(1+o(1))+r^2n\vartheta_n^{2r}(\kappa_n-x_n)^2}$$
$$\leqslant \frac{8(1+o(1))}{8(1+o(1))+r^2(1-2r)nx_n^{(r+q)/q}}.$$
Since the assumption on $x_n$ implies $nx_n^{(r+q)/q}\to\infty$ this implies $p_n\to 0$.

By Mogulskii's inequality and the above we get
$$P_0^n(V_n\geqslant \sqrt{n}x_n)(1-e^{-nx_n^2})\geqslant \exp\{-nx_n^{(r+q)/q}\log(x_n^{(r+q)/q}/\vartheta_n)(1+o(1))-nx_n^2p_n\}-e^{-nx_n^2}.$$
Observe that $nx_n^{(r+q)/q}\log(x_n^{(r+q)/q}/\vartheta_n)/nx_n^2=x_n^{(r-q)/q}\log(x_n^{(r+q)/q}/\vartheta_n)\to 0$ by the assumption on $x_n$. Therefore the second term on the right hand side of the last estimate is of higher order than the first. Logarithming both sides and dividing by $nx_n^2$ gives (7). \hfill $\Box$\\

\noindent{\bf Appendix. Classical moderate deviation theorem}\\

In this section we reprove the classical MD theorem for i.i.d. random variables using Mogulskii's inequality. We do this to evidence strong similarity of the proofs of Theorems 2 and 4.

Let $\xi_1,\xi_2,...$ be a sequence of i.i.d. real random variables with distribution $P$, $E\xi_1=0,$ Var\,$\xi_1=1$ and $\varphi(\lambda)=Ee^{\lambda \xi_1}$ finite for $\lambda\in [0,\Lambda],\;\Lambda>0$.\\

{\bf Theorem 4.} {\it If $x_n\to 0$ is such that $nx_n^2\to\infty$ then we have}
$$ -\lim_{n\to\infty} \frac{1}{nx_n^2}\log Pr\left(\frac{1}{\sqrt{n}}\sum_{i=1}^n \xi_i\geqslant \sqrt{n}x_n\right)=\frac{1}{2}.$$

{\bf Proof.}

{\it Upper estimate}. The function $\varphi(\lambda)$ is analytic on $[0,\Lambda]$ and can be written in a form 
$$ \varphi(\lambda)=1+\frac{\lambda^2}{2}\psi(\lambda),$$
where $\psi(\lambda)$ is analytic, $\psi(\lambda)\geqslant 0$ and $\psi(0)=1$. By independence and Markov's inequality we get for arbitrary $\lambda\in [0,\Lambda]$
$$Pr\left(\frac{1}{\sqrt{n}}\sum_{i=1}^n \xi_i\geqslant \sqrt{n}x_n\right)=Pr\left(\prod_{i=1}^n e^{\lambda \xi_i}\geqslant e^{n\lambda x_n}\right)\leqslant e^{-n\lambda x_n}\varphi^n(\lambda).$$
Setting $\lambda=x_n$, logarithming and dividing by $nx_n^2$ we obtain from the form of $\varphi(\lambda)$
$$\frac{1}{nx_n^2}\log Pr\left(\frac{1}{\sqrt{n}}\sum_{i=1}^n \xi_i\geqslant \sqrt{n}x_n\right)\leqslant -1+\frac{\log(1+\frac{x_n^2}{2}\psi(x_n))}{x_n^2}$$
which immediately implies
$$\limsup_{n\to\infty}\frac{1}{nx_n^2}\log Pr\left(\frac{1}{\sqrt{n}}\sum_{i=1}^n \xi_i\geqslant \sqrt{n}x_n\right)\leqslant -\frac{1}{2}.$$

{\it Lower estimate.} For any $\lambda\in[0,\Lambda]$ consider the distribution $Q_{\lambda}\ll P$ defined by $\frac{dQ_{\lambda}}{dP}(y)=e^{\lambda y}/\varphi(\lambda)$. Then
$$m(\lambda)=\int y dQ_{\lambda}=\frac{1}{\varphi(\lambda)}\int ye^{\lambda y}dP(y)=\frac{\varphi'(\lambda)}{\varphi(\lambda)}$$
and the Kullback-Leibler distance of $Q_{\lambda}$ from $P$ can be expressed by
$$D(Q_{\lambda}||P)=\int\frac{1}{\varphi(\lambda)}e^{\lambda y}(\lambda y-\log \varphi(\lambda))dP(y)=\lambda\frac{\varphi'(\lambda)}{\varphi(\lambda)}-\log\varphi(\lambda).$$
For $n\geqslant 1$ and $\varepsilon\in (0,1/3)$ let $\lambda_n>0$ be such that $m(\lambda_n)=(1+\varepsilon)x_n$.
Since $\log \varphi(\lambda)$ is strictly convex then the function $m(\lambda)=\varphi'(\lambda)/\varphi(\lambda)$ is increasing and $m(0)=0$. Hence $\lambda_n\to 0$. For $n$ sufficiently large i.e. such that $|\psi(\lambda_n)-1|<\varepsilon/8,$ $|\lambda_n\psi'(\lambda_n)|<\varepsilon/4$ and $\lambda_n^2\psi(\lambda_n)<\varepsilon/4$ we have
$$ (1+\varepsilon)x_n=m(\lambda_n)=\frac{\varphi'(\lambda_n)}{\varphi(\lambda_n)}=\frac{\lambda_n\psi(\lambda_n)+\lambda_n^2\psi'(\lambda_n)/2}{1+\lambda_n^2\psi(\lambda_n)/2}\leqslant \lambda_n(1+\varepsilon/4)\leq \lambda_n(1+\varepsilon)$$
and similarly
$$ (1+\varepsilon)x_n\geqslant \lambda_n\frac{1-\varepsilon/4}{1+\varepsilon/8}$$
which implies 
$$x_n\leqslant\lambda_n\leqslant x_n \frac{(1+\varepsilon)(1+\varepsilon/8)}{1-\varepsilon/4}\leqslant (1+2\varepsilon)x_n.\eqno (21)$$

For $\lambda_n$ defined above we have
$$ D(Q_{\lambda_n}||P)=(1+\varepsilon)\lambda_nx_n-\log\varphi(\lambda_n).$$
In Mogulskii's inequality (Theorem A) set $Q=Q_{\lambda_n},\;A=[x_n,\infty),\; M=2nx_n^2$.  Since $\varphi''(0)=1$ then the variance of $Q_{\lambda_n}$ is equal to $\rho_n^2=\varphi''(\lambda_n)/\varphi(\lambda_n)-m^2(\lambda_n)\to 1$. Hence for $n$ sufficiently large, by the assumption $nx_n^2\to\infty$ and from Cantelli's inequality we obtain
$$ p_n=Pr(\eta_1+...+\eta_n< nx_n)=Pr\left(\sum_{i=1}^n(\eta_i-m(\lambda_n))< -\varepsilon nx_n\right)\leqslant \frac{n\rho_n^2}{n\rho_n^2+\varepsilon^2n^2x_n^2}\to 0.$$
From (21) we have $\lambda_n^2\geq \lambda_nx_n$ and for $n$ sufficiently large $\log(1+\lambda_nx_n\psi(\lambda_n)/2)\geq (1-\varepsilon)\lambda_nx_n/2$. Hence, again (21) and Mogulskii's inequality imply
$$ Pr(\xi_1+...+\xi_n\geqslant nx_n)(1-e^{-2nx_n^2})$$
$$\geqslant \exp\{-(1+\varepsilon)n\lambda_nx_n+n\log(1+\lambda_n^2\psi(\lambda_n)/2)-2nx_n^2p_n\}-e^{-2nx_n^2}$$
$$\geqslant \exp\left\{-\frac{1+3\varepsilon}{2}n\lambda_nx_n-2nx_n^2p_n\right\}-e^{-2nx_n^2}\geqslant \exp\left\{(-\frac{1}{2}-\frac{7}{2}\varepsilon)nx_n^2-2nx_n^2p_n\right\}-e^{-2nx_n^2}.$$
Logarithming and dividing by $nx_n^2$ both sides we obtain
$$\frac{1}{nx_n^2}\log Pr(\xi_1+...+\xi_n\geqslant nx_n)\geqslant-\frac{1}{2}-\frac{7}{2}\varepsilon +o(1)$$
which, due to arbitrariness of $\varepsilon\in (0,1/3)$, gives
$$\liminf_{n\to\infty}\frac{1}{nx_n^2}\log Pr\left(\frac{1}{\sqrt{n}}\sum_{i=1}^n\xi_i\geqslant \sqrt{n}x_n\right)\geqslant -\frac{1}{2}$$
and finishes the proof. \hfill $\Box$

\vspace{8mm}
\noindent Tadeusz Inglot\\
Faculty of Pure and Applied Mathematics,\\
Wroc{\l}aw University of Science and Technology,\\
Wybrze\.ze Wyspia\'nskiego 27, 50-370 Wroc{\l}aw, Poland.\\
E-mail: Tadeusz.Inglot@pwr.edu.pl

\end{document}